\documentclass[10pt]{amsart}

\usepackage{amsmath, amsfonts, amssymb, amsthm, graphics}

\newtheorem{theorem}{Theorem}[section]

\newtheorem{proposition}[theorem]{Proposition}

\theoremstyle{definition}
\newtheorem{definition}[theorem]{Definition}

\numberwithin{equation}{section}

\begin{document}

\title{Hardy and Littlewood theorems and the Bergman distance}

\author{Marijan Markovi\'{c}}

\dedicatory{Dedicated to the memory of Professor Miroslav Pavlovi\'{c}}

\begin{abstract}
We obtain non-Euclidean versions of classical theorems due to Hardy and Littlewood concerning smoothness of the boundary function of an analytic
mapping on the unit disk with an appropriate growth condition.

\medskip

\noindent\textsc{R\'{e}sum\'{e}.}
Nous obtenons des versions non euclidiennes des th\'{e}or\`{e}mes classiques dus \`{a} Hardy et Littlewood concernant la r\'{e}gularit\'{e} de
la  fonction fronti\`{e}re d'une fonction analytique sur le  disque unit\'{e} avec une condition de croissance appropri\'{e}e.
\end{abstract}

\address{Faculty of Sciences and Mathematics\endgraf University of Montenegro\endgraf
D\v{z}ord\v{z}a Va\v{s}ingtona bb\endgraf
81000  Podgorica, Montenegro}

\email{marijanmmarkovic@gmail.com}



\maketitle

\section{Introduction}
For the theory of Hardy spaces $H^p$ on the unit disk $\mathbb{D}$ we refer to the classical Duren book \cite{DUREN.BOOK} and the recent Pavlovi\'{c}
monograph \cite{PAVLOVIC.BOOK}. Recall that the space $H^\infty$ contains bounded analytic mappings on $\mathbb{D}$. The space $H^p$, $p\in [1,\infty)$,
consists of analytic mappings $f$ on $\mathbb{D}$ such that  $m_p(r,|f|)$  remains bounded in $r\in (0,1)$, where for a continuous and nonnegative
function $g$ on the unit disk we have  denoted
\begin{equation*}
m_p(r,g) = \left(\int_{-\pi}^{\pi}  g(r e^{it})^p dt\right)^{\frac 1p}.
\end{equation*}

It is well known that an analytic mapping  $f\in H^p$, $p\in [1,\infty]$,  has radial boundary value
\begin{equation*}
f_b(t) = f(e^{it})=  \lim_{r\to 1} f(r e^{it})\  \text{for a.e.}\  t\in [-\pi,\pi].
\end{equation*}
Moreover, the boundary function $f_b$ belongs to the Lebesgue space $L^p [-\pi ,\pi]$,  and the mapping $f$ is the Poisson extension of  $f_b$, i.e.,
\begin{equation*}
f(z)=\mathrm {P}[f_b (t)] (z)=\frac 1{2\pi}\int_ {-\pi}^\pi \frac {1-r^2}{1-2r\cos (\theta - t) + r^2} f_b (t) dt,\quad z = re^{i\theta}\in \mathbb{D}.
\end{equation*}

The Lipschitz class $\Lambda_\alpha$, $\alpha \in (0,1]$, consists of $2\pi$-periodic functions $\phi $ on  $\mathbb{R}$ such that
\begin{equation*}
\sup _{|s-t|<h}| \phi (t) - \phi(s)| = \mathcal {O}(h^\alpha),\quad h\to 0.
\end{equation*}

The mean Lipschitz class $\Lambda _\alpha^p$, $\alpha \in (0,1]$, $p\in [1,\infty)$, contains $2\pi$-periodic functions  $\phi$ defined a.e. on
$\mathbb{R}$ such that
\begin{equation*}
\sup _{s\in (0,h)}\left(\int_{-\pi} ^{\pi}|\phi (t+s) - \phi (t)|^p dt\right)^{\frac 1p} = \mathcal {O}(h^\alpha),\quad h\to 0.
\end{equation*}

The following two propositions are classical results due to Hardy and Littlewood \cite{HL.MZ}. We refer to the fifth chapter of  \cite{DUREN.BOOK}.

\begin{proposition}\label{PROP.HL1}
An analytic mapping $f$ on $\mathbb{D}$ satisfies $|f'|(z)=\mathcal{O}(1-|z|)^{\alpha-1}$, $|z|\to 1$, $\alpha\in (0,1]$, if and only if it has
continuous  extension on  $\overline{\mathbb{D}}$ and $f_b\in \Lambda_\alpha$.
\end{proposition}

\begin{proposition}\label{PROP.HL2}
An analytic mapping $f$ on $\mathbb{D}$ satisfies $m_p(r,|f'|)=\mathcal {O}(1-r)^{\alpha-1}$, $r\to 1$, $\alpha\in (0,1]$, $p\in [1,\infty)$, if
and only if  $f\in H^p$ and $f_b\in\Lambda^p_\alpha$.
\end{proposition}

In 1982, Yamashita gave the hyperbolic counterparts of Proposition \ref{PROP.HL1} and Proposition  \ref{PROP.HL2};
see \cite[Theorem 1 and Theorem 2]{YAMASHITA.TAMS}.  Instead  of the growth of  $|f'|$ his theorems  concern  the growth of the hyperbolic
derivative
\begin{equation*}
f^\ast (z) = \frac { |f'(z)| }{1-|f(z)|^2},\quad z\in \mathbb{D},
\end{equation*}
of an analytic mapping $f:\mathbb{D}\to \mathbb{D}$. On the other hand, instead of the standard Lipschitz classes,  the  Yamashita results  contain
their hyperbolic  analogues.

The aim of this article is to consider the Yamashita approach in a more general than the hyperbolic setting on the unit disk. We consider analytic
mappings of $\mathbb{D}$ into a domain in $\mathbb{C}$ equipped with a metric with ceratin  properties. Particularly, this domain may be bounded and
equipped with  the  Bergman metric.

\section{Preliminaries}
We say that a function $\rho$ is a metric on a domain $D\subseteq\mathbb{C}$ if it is positive and continuous on $D$. The $\rho$-distance on $D$ is
\begin{equation*}
d(\zeta, \eta) = \inf_\gamma \int _\gamma \rho(\zeta) |d\zeta| = \int_a^b \rho (\gamma(t)) |\gamma'(t)| dt,\quad \zeta,\eta\in  D,
\end{equation*}
where the infimum is taken over all  piecewise  $\mathcal{C}^1$-smooth curves  $\gamma:[a,b]\to D$ with endpoints at $\zeta$  and $\eta$, i.e.,
$\gamma(a)=\zeta$  and  $\gamma(b)=\eta$.

It may  be shown that
\begin{equation*}
\lim_{\eta \to  \zeta } \frac {d  (\zeta,\eta)}{|\zeta - \eta|} = \rho (\zeta),\quad \zeta\in D.
\end{equation*}
For this result we refer to \cite{MARKOVIC.JGA,ZHU.JLMS}.

For the notion of Bergman metric  we  refer to the 4th  chapter in  the  Krantz book \cite{KRANTZ.BOOK}, but for  the sake of completeness we mention
below some facts.

For a bounded  domain $D$ let
\begin{equation*}
A^2(D) = \{f: f \  \text{is analytic mapping on $D$ and   $\int_D |f|^2 dA$ is finite}\},
\end{equation*}
be the  Bergman space; $dA$ is  the area measure. In  $A^2 (D)$ one may introduce  the inner product
\begin{equation*}
\left<f,g\right> _{A^2(D)} = \int _\Omega f(z) \overline {g(z)}dA (z),\quad f,g\in A^2(D).
\end{equation*}
With this inner product  $A^2 (D)$ is the Hilbert space.

Among other approaches, one can use abstract Hilbert space methods to prove the existence of the so called Bergman kernel
$K_D: D\times D \to \mathbb {C}$  which has  the following  properties:

(1) $K(z,w) = \overline {K (w,z)}$,  $(z, w) \in D\times D$;

(2) $K_D^z (w) = K(z ,w) \in  A^2(D)$ for every $z\in D$;

(3) the  reproducing property: for every $f \in A^2(D)$ holds
\begin{equation*}
f(z) = \left <K_D^z,f\right> _{A^2(D)} = \int_D K_D(z,w) f(w) dA(z),\quad z\in D.
\end{equation*}

The Hilbert space $A^2 (D)$ is separable, so  it has a countable  base  $\{\phi_j:j \in \mathbb{N} \}$. The Bergman kernel  may be represented  as
\begin{equation*}
K (z,w) = \sum _{j =1 }^\infty \phi_j (z) \overline {\phi_j (w)},\quad z,w\in D.
\end{equation*}
This  possibility  gives opportunity to show that the continuous  function
\begin{equation*}
\rho_D (\zeta)=\sqrt  {\frac {\partial ^2}{\partial \zeta \partial\overline {\zeta}} \log  K_D(\zeta,\zeta)},\quad  \zeta\in D.
\end{equation*}
is positive on  $D$. Therefore,  it is a metric on  $D$   which is called  the Bergman metric.

It is well known that  the  Bergman metric on the unit disk coincides up to a multiplicative  constant  with the hyperbolic metric on  $\mathbb{D}$
which is  given by
\begin{equation*}
\rho (\zeta) = \frac 1{ 1-|\zeta|^2},\quad \zeta\in \mathbb{D}.
\end{equation*}
The  hyperbolic distance on $\mathbb{D}$  is  explicitly given  by
\begin{equation*}
\sigma(\zeta,\eta) = \frac 12 \log \frac{|1-\overline{\zeta} \eta| +|\zeta - \eta|}
{|1-\overline{\zeta}\eta| -|\zeta - \eta|},\quad \zeta,\eta\in\mathbb{D}.
\end{equation*}

Introduce now the generalized Lipschitz classes -- non-Euclidean analogues to $\Lambda_\alpha$ and $\Lambda_\alpha^p$. Let $\rho$ be a metric on a
domain $D$ and let $d = d_\rho$ be the $\rho$-distance on the same domain. We denote by $\rho\Lambda_\alpha$, $\alpha\in(0,1]$, the class of all
$2\pi$-periodic functions $\phi$ on $\mathbb{R}$ with values in $D$ such that
\begin{equation*}
\sup_{|s-t|<h} d ( \phi (t),\phi (s)) = \mathcal {O}(h^\alpha),\quad h\to 0.
\end{equation*}
By $\rho\Lambda_\alpha^p$, $\alpha\in (0,1]$, $p\in [1,\infty)$, we denote the class of $2\pi$-periodic function defined a.e. on $\mathbb{R}$ with
values in the domain $D$ such that
\begin{equation*}
\sup_{s\in (0,h)} \left(\int_0 ^{2\pi}d  (\phi(t+s),\phi(t))^p {dt} \right)^{\frac 1p} = \mathcal {O}(h^\alpha),\quad h\to 0.
\end{equation*}
Note that, if $p<p'$, then $\rho\Lambda_\alpha\subseteq \rho\Lambda _\alpha^{p'}\subseteq \rho\Lambda _\alpha^p$, $\alpha\in (0,1]$, $p,p'\in [1,\infty)$.
If $\rho$ is the hyperbolic metric on $\mathbb{D}$,  we have the hyperbolic Lipschitz classes  $\sigma\Lambda_\alpha$ and $\sigma\Lambda_\alpha^p$.

Let $\Omega$ be a domain in $\mathbb{C}$ or in $\mathbb{R}^1$. The $\rho$-derivative of a mapping  $f:\Omega\to D$ at $z\in\Omega$ is defined by
\begin{equation*}
f^\ast  (z) = \lim _{w\to z}\frac {d(f(z), f(w))}{|z-w|}.
\end{equation*}
If $f:\Omega\to D$ is analytic (differentiable) mapping at $z_o\in\Omega$, then the $\rho$-derivative of $f$ at $z_o$  may be expressed  in the
following  way
\begin{equation*}
f^\ast(z_o) = \rho (f(z _o)) |f'(z_o)|.
\end{equation*}
Indeed,   if $f' (z_o)\ne 0$, then  $f(z)\ne f(z_o)$, if $z$ is in a  sufficiently small neighborhood of $z_o$; namely, if $f(z_n)= f(z_o)$ for a
sequence  $z_n\to z_o$, then
\begin{equation*}
f'(z_o)  = \lim _{n\to \infty} \frac{ f(z_n) - f(z_o) }  { z_n - z_o }  = 0.
\end{equation*}
It follows
\begin{equation*}\begin{split}
f^\ast (z_o) & = \lim_{z\to z_o} \frac { d (f(z),f(z_o)) }{|z-z_o|}
= \lim_{z\to z_o} \frac {d(f(z),f(z_o))}{|f(z)-f(z_o)|} \lim_{z\to z_o} \frac{|f(z)-f(z_o)|}{|z-z_o|}
\\& =\rho(f(z_o)) | f'(z_o) |.
\end{split}\end{equation*}
Now, let $f'(z_o) = 0$,  and let $r$ be a positive  number such that the closed disk $\overline{\mathbb{D}}( f(z_o),r)$ is contained in $D$. Denote
$m=\max_{\zeta\in\overline{\mathbb{D}}( f(z_o),r)} \rho (\zeta)$. For  $z$ sufficiently close to $z_o$ such  that $|f(z) - f(z_o)|\le r$  we   have
\begin{equation*}
d (f(z),f(z_o))\le m |f(z)  -  f(z_o)|.
\end{equation*}
It follows
\begin{equation*}
f^\ast (z_o)
= \lim _{z \to z_o } \frac{ d  (f(z),f(z_o))}{|z-z_o|}\le m \lim _{z \to z_o }\frac{ |f(z)- f(z_o) | }{|z-z_o|} = 0,
\end{equation*}
which proves this case.

Let us note that for an analytic mapping $f:\mathbb{D}\to\mathbb{D}$, and for the hyperbolic metric on $\mathbb{D}$, $f^ \ast(z)$ is the hyperbolic
derivative of  $f$ at $z\in \mathbb{D}$.

\section{Two  generalizations   of the  Yamashita result}
We  will state  the  Yamashita results \cite[Theorem 1 and Theorem 2]{YAMASHITA.TAMS} in  the following  proposition.

\begin{proposition}\label{PROP.YAMASHITA}
Let $f:\mathbb{D}\to \mathbb{D}$ be an analytic mapping,                        and let $f^\ast$ be the hyperbolic derivative of $f$. Then  we have
\begin{equation*}
f^\ast(z) = \mathcal {O}(1-|z|)^{\alpha-1}, |z|\to 1,\quad   \alpha\in (0,1],
\end{equation*}
if  and only if $f_b\in \sigma\Lambda_\alpha$, and
\begin{equation*}
m_p(r , f^\ast ) = \mathcal {O}(1-r)^{\alpha-1},  r\to 1,\quad   \alpha\in (0,1],  p\in [1,\infty),
\end{equation*}
if and only if $f_b\in \sigma\Lambda^p_\alpha$.
\end{proposition}

Our  main results are given in the  fourth and fifth  section. Here we will state their consequences for bounded domains with the Bergman metric and
for Dini-smooth Jordan domains with the quasi-hyperbolic metric which both extend Proposition \ref{PROP.YAMASHITA}. The Bergman metric on the unit
disk $\mathbb{D}$  coincides up to a positive constant with the  hyperbolic metric on this domain. Therefore,  the result  given below  indeed
generalize the Yamashita one.

As a direct consequence of  Theorem  \ref{TH.1}, Theorem \ref{TH.2}, and Theorem \ref{TH.3}  we  have the following theorem.

\begin{theorem}\label{TH.BERGMAN}
Let $D\subseteq \mathbb{C}$  be   a bounded  domain equipped with the Bergman metric    $\rho_D$. For an  analytic mapping  $f:\mathbb{D}\to D$ let
\begin{equation*}
f^\ast (z) = \rho_D (f(z)) {| f'(z) | },\quad    z\in  \mathbb {D},
\end{equation*}
be the $\rho_D$-derivative of $f$. Then
\begin{equation*}
f^\ast (z) = \mathcal {O}(1-|z|)^{\alpha-1}, |z|\to 1,\quad \alpha\in (0,1],
\end{equation*}
if and only if $f_b\in \rho_D\Lambda_\alpha$, and
\begin{equation*}
m_p(r , f^\ast ) = \mathcal {O}(1-r)^{\alpha-1},  r\to  1, \quad \alpha\in (0,1],p\in [1,\infty),
\end{equation*}
if and only if $f_b  \in \rho_D\Lambda^p_\alpha$.
\end{theorem}

A $\mathcal {C}^1$-smooth Jordan curve $\gamma$is called  Dini-smooth if, considered as a function of the arc length over the segment $[0,l]$, where
$l$ is the length of  $\gamma$,  satisfies
\begin{equation*}
| \gamma' (s ) - \gamma' (t)|\le \omega(|s - t|),\quad  s,t\in [0,l],
\end{equation*}
and  $\omega(x)$ is an increasing function  such that
\begin{equation*}
\int_0^l \frac { \omega(s)}s
\end{equation*}
is finite.

\begin{theorem}\label{TH.JORDAN}
Let $D$  be   a  Jordan domain with Dini-smooth boundary  with the quasi-hyperbolic metric
\begin{equation*}
\rho(\zeta) = \frac 1 {d(\zeta)},\quad  \zeta\in D,
\end{equation*}
where  $d(\zeta)$ is the Euclidean  distance of $\zeta$ form $\partial D$. For an analytic mapping $f:\mathbb{D}\to D$ let
\begin{equation*}
f^\ast_q(z) = \frac { |f'(z)|} { d(z)},\quad z\in \mathbb{D}
\end{equation*}
be  the quasi-hyperbolic derivative of $f$. We have
\begin{equation*}
f^\ast_q (z) = \mathcal {O}(1-|z|)^{\alpha-1}, |z|\to 1,\quad \alpha\in (0,1],
\end{equation*}
if and only if $f_b\in \rho\Lambda_\alpha$, and
\begin{equation*}
m_p(r , f^\ast_q ) = \mathcal {O}(1-r)^{\alpha-1},  r\to  1, \quad \alpha\in (0,1],p\in [1,\infty),
\end{equation*}
if and only if $f_b  \in \rho\Lambda^p_\alpha$.
\end{theorem}

\begin{proof}
Let the  quasi-hyperbolic  distance on the  domain $D\subseteq \mathbb{C}$  be denoted $d_q$, and denote by  $d_\beta$ the Bergman distance  on the
same domain. It is possible to estimate  the Bergman metric via the quasi-hyperbolic  metric on a Jordan domain with   Dini-smooth. Quite  recently
Nikolov and Trybula \cite{NIKOLOV.COLLM} proved: There exists a  constant $c>1$  such  that
\begin{equation*}
\sqrt{2}\log\left(1+\frac{|\zeta-\eta|} {c\sqrt{d(\zeta)d(\eta)}} \right)\le d_\beta  (\zeta,\eta)\le
\sqrt{2}\log\left(1+\frac{c|\zeta-\eta|} { \sqrt{d(\zeta)d(\eta)}} \right), \quad  \zeta, \eta\in D.
\end{equation*}
This forces the following relation between the Bergman and the quasi--hyperbolic metric (if we divide each side above by  $|\zeta - \eta|$ and  let
$\eta\to \zeta$):
\begin{equation}\label{EQ.QH}
\frac {c^{-1}\sqrt{2}} {d(\zeta )}\le   \rho_D(\zeta)\le   \frac {c\sqrt{2}}{d(\zeta )},\quad \zeta\in D.
\end{equation}
It follows
\begin{equation*}
{c}^{-1}{\sqrt{2}} d_{q} (\zeta,\eta) \le  d_\beta  (\zeta,\eta) \le c\sqrt{2}   d_{q} (\zeta,\eta),\quad \zeta,  \eta\in D.
\end{equation*}
Therefore, we have the coincidence of the classes  $\rho_D \Lambda_\alpha = \rho  \Lambda_\alpha$, $\alpha\in (0,1]$, and
$\rho_D \Lambda_\alpha ^p=\rho\Lambda_\alpha^p$, $\alpha\in (0,1]$, $p\in [1,\infty)$.
Clearly, the  inequality \eqref{EQ.QH} implies the relation between the derivatives:
\begin{equation*}
c^{-1}\sqrt{2} f_q^\ast(z)\le f^\ast (z)\le c\sqrt{2} f_q^\ast(z),\quad z\in \mathbb {D},
\end{equation*}
where $f^\ast$ is the $\rho_D$-derivative.   The conclusion of this theorem now follows from   Theorem \ref{TH.BERGMAN}.
\end{proof}

\section{Growth of the derivative and the  boundary  function}
In the theorems which follows we assume that $D$ is a domain equipped with a metric $\rho$ and the corresponding distance $d = d_\rho$,  such that
the following two  conditions are satisfied:

(1) $C_\rho : = \inf_{\zeta \in D} \rho (\zeta)  >0$;

(2) $d (\zeta, \eta) \to \infty $,  if $\zeta\to\partial D$;  $\eta\in D$ is arbitrary.

Note that (1) is satisfied if the metric $\rho$ has the stronger property: $\rho(\zeta)\to\infty$, if $\zeta\to\partial D$ (the boundary is in the
topology of $\overline{\mathbb{C}}$). For example this stronger condition is satisfied by the Bergman metric on a bounded domain. Also, it is known
that the  Bergman metric   has the second  property (see \cite{KRANTZ.BOOK}).

The first condition  immediately  implies that
\begin{equation}\label{RO.E}
d(\zeta,\eta)\ge  C_\rho |\zeta - \eta|,\quad  \zeta, \eta\in D,
\end{equation}
which forces the inclusions $\rho\Lambda_\alpha\subseteq\Lambda_\alpha$, $\rho\Lambda^p_\alpha\subseteq\Lambda^p_\alpha$, $\alpha\in (0,1]$,
$p\in [1,\infty)$.

The $\rho$-Hardy space $\rho H^p$, $p\in [1,\infty)$, is defined as  a  space of  analytic mappings $f:\mathbb{D}\to D$ such that
\begin{equation*}
\|f\|_{\rho H^p} = \sup_{r\in (0,1)} \left(\int_{-\pi}^\pi d (f(re^{it}),f(0))^pdt\right)^\frac 1p
\end{equation*}
is finite. The property \eqref{RO.E} implies the inclusion $\rho H^p \subseteq  H^p$, $p\in [1,\infty)$. Indeed, for $f\in\rho H^p$ and $r\in(0,1)$
we have
\begin{equation*}
\left(\int _{-\pi}^\pi  |f(re^{it}) -  f(0) |^p dt\right)^{\frac 1p} \le {C_\rho}^{-1} {\|f\|_{\rho H^p}},
\end{equation*}
which throughout some elementary inequalities    easily   implies that $f\in H^p$, since
\begin{equation*}
m_p (|f|,r)\le 4\pi |f(0)| + 2{C_\rho}^{-1}   {\|f\|_{\rho H^p}}.
\end{equation*}

The  main results of this section  are stated in the following  two  theorems. Note that if  we take  $\rho\equiv 1$  and $D = \mathbb{C}$
(then $d_\rho$ is   the Euclidean distance) Theorem \ref{TH.1} became Proposition \ref{PROP.HL1}, and  Theorem \ref{TH.2}  became  one part of
Proposition \ref{PROP.HL2}.

\begin{theorem}\label{TH.1}
An analytic mapping $f:\mathbb{D}\to D$ satisfies the condition
\begin{equation*}
f^\ast(z) = \mathcal {O}(1-|z|)^{\alpha-1},|z|\to 1,\quad \alpha\in(0,1].
\end{equation*}
if and only if it has continuous extension on $\overline{\mathbb{D}}$ and $f_b\in\rho\Lambda_\alpha$.
\end{theorem}

\begin{theorem}\label{TH.2}
If an analytic mapping  $f:\mathbb{D}\to D$  satisfies
\begin{equation}\label{EQ.GROWTH2}
m_p(r, f^\ast) = \mathcal {O}(1-r)^{\alpha-1}, r\to 1,\quad\alpha\in(0,1], p\in[1,\infty)
\end{equation}
then $f\in \rho H^p$  and  $f_b\in\rho\Lambda^p_\alpha$.
\end{theorem}

The following fact will be used several times in the  proofs of these theorems: For an analytic  mapping  $f:\mathbb {D}\to D$ we  have
\begin{equation*}
d  (f(z),f(w))\le\int_\gamma f^\ast  (\zeta)|d\zeta|,
\end{equation*}
where $\gamma:[a,b]\to\mathbb{D}$ is a piecewise $\mathcal{C}^1$-smooth curve which joins  $z$ and $w$. Indeed, since $f\circ\gamma:[a,b]\to D$  is
a  piecewise      $\mathcal{C}^1$-smooth curve   which  connects  $f(z)$ and $f(w)$,   we have
\begin{equation*}\begin{split}
d (f(z),f(w)) &  \le\int_{f\circ \gamma} \rho (\eta)|d\eta| = \int_a^b \rho ( f (\gamma (t)) )| f'(\gamma (t))| |\gamma'(t)|dt
\\& =   \int_{ \gamma} \rho (f(\zeta))|f'(\zeta)||d\zeta | = \int_\gamma f^\ast  (\zeta)|d\zeta|,
\end{split}\end{equation*}
which we aimed to prove.

\begin{proof}[Proof of Theorem \ref{TH.1}]
Let  us first assume that $f^\ast (z) = \mathcal {O}  (1- |z|)^{\alpha-1}$, $|z|\to 1$. Since
\begin{equation*}
f^\ast (z) = \rho ( f(z) ) |f'(z)|\ge C_\rho |f' (z)|,\quad   z\in \mathbb{D},
\end{equation*}
we also have $|f'(z)| = \mathcal{O}(1- |z|)^{\alpha-1}$, $|z|\to 1$. By Proposition \ref{PROP.HL1}, we conclude that $f$ has continuous extension
on $\overline{\mathbb{D}}$, and moreover $f_b\in\Lambda_\alpha$. It remains to show that $f_b\in \rho\Lambda_\alpha$.

Let $C$ be  a constant  such that  $f^\ast(z)\le C(1- |z|)^{\alpha-1} $, $z\in \mathbb{D}$. For  $r\in (0,1)$ and $t\in [-\pi,\pi]$ we have
\begin{equation*}\begin{split}
d(f(re^{it}), f(0)) & \le \int_0^r f^\ast (x  e^{it}) dx \le \int_0^1 f^\ast (x  e^{it})dx
\\& \le C  \int _0^1 (1-x)^{\alpha-1} dx =\frac {C}\alpha.
\end{split}\end{equation*}
Having in mind  our assumption on the distance $d$, it follows that
\begin{equation*}
f_b(t) = \lim_ {r\to 1} f(re^{it})\in D, \quad t\in [-\pi,\pi].
\end{equation*}

Let $\tilde d$ be the half of the distance between $\{f_b(t):t\in [-\pi,\pi]\}$ and $\partial D$ (in the case $D=\mathbb{C}$, take $\tilde{d}=1$).
Denote
\begin{equation*}
M = \max_{\zeta\in \bigcup_{t\in [-\pi,\pi]} \overline{\mathbb{D}}(f_b(t),\tilde {d})} \rho ( \zeta).
\end{equation*}
Let $h>0$ be such that  $Ch^\alpha<\tilde{d}$. Because of uniform continuity of $f_b$, for sufficiently small $h$  we have:
$|t-s|<h \Rightarrow |f_b (t)-f_b (s)|\le Ch^\alpha<\tilde {d}$.
It follows
\begin{equation*}\begin{split}
d  (f_b(t), f_b(s))& = \inf_\gamma \int_\gamma \rho (\zeta) |d\zeta|
\le  \int_{[f_b(t),f_b(s)]}\rho (\zeta) |d\zeta| \\&\le M  |f_b (t) -  f_b  (s)|\le MC h^\alpha
\end{split}\end{equation*}
(let us say that $\gamma\subseteq D$ is among all piecewise $\mathcal{C}^1$-smooth curves that join $f_b(t)$ and $f_b(s)$, and $[f_b(t),f_b(s)]$
is the  segment  $\{f_b  (t)  + \lambda (f_b (s)- f_b(t)), \lambda \in [0,1]\}$). This proves that $f_b\in \rho \Lambda_\alpha$.

In order to prove the converse, assume that $f$ has continuous extension on $\overline{\mathbb{D}}$ with values in $D$, and
$f_b\in\rho\Lambda_\alpha$. Then $f_b\in\Lambda_\alpha$, so we may apply the Hardy and Littlewood result given in Proposition \ref{PROP.HL1}.
We conclude that there exists a constant $C$ such that  $|f'(z)|\le  C (1- |z|)^{\alpha-1}$, $z\in\mathbb {D}$. Since $f_b(t) =  f(e^{it})$ is
continuous, let $m = \max_{t\in [0,2\pi]}\rho( f_b(t) )$. By uniform continuity of $\rho \circ f$ on $\overline{\mathbb{D}}$, if $r$ is
sufficiently close to  $1$, we have $\rho (f(z))\le 2m$, $z\in \overline{\mathbb {D}} \backslash r\mathbb {D}$.  It follows that
\begin{equation*}
f^\ast (z) = \rho (f(z)) |f'(z)|\le 2 m C (1- |z|)^{\alpha-1},\quad   |z|>r,
\end{equation*}
i.e., $f^\ast (z) = \mathcal {O}(1-|z|)^{\alpha-1}$, $|z|\to 1$, which we aimed to prove.
\end{proof}

\begin{proof}[Proof of Theorem \ref{TH.2}]
We shall first show that $f\in \rho H^p$ which implies that $f_b(t)$ is finite for a.e. $t\in [-\pi,\pi]$. Then we prove that  $f_b(t)\in D$ for a.e.
$t\in [-\pi,\pi]$.

Assume that for a constant $C$ there holds
\begin{equation}\label{EQ.CALPHA}
m_p (r,f^\ast)\le C (1-r)^{\alpha-1},\quad  r\in (0,1).
\end{equation}
Since
\begin{equation*}
d  (f(re^{it}), f(0)) \le  \int_0 ^r f ^\ast (x e^{it})dx,
\end{equation*}
applying  the   Minkowski inequality, we obtain
\begin{equation*}\begin{split}
\left(\int _{-\pi}^\pi d (f(re^{it}), f(0))^p dt\right)^{\frac 1p} &
\le \left( \int _{-\pi}^\pi    \left(\int_0 ^r f  ^\ast (x e^{it})dx\right)^pdt\right)^{\frac 1p}
\\&\le\int_0 ^r \left( \int _{-\pi}^\pi   f  ^\ast (x e^{it})^p dt\right)^{\frac 1p}dx
\\&\le\int_0 ^r m_p (x, f ^\ast)dx
\\& \le C \int_0 ^1 (1-x)^{\alpha-1} dx
\\&=\frac C \alpha.
\end{split}\end{equation*}
Therefore,  $f\in \rho H^p$. By the Fatou theorem and the above inequality we conclude that
\begin{equation*}
\int _{-\pi}^\pi\liminf_{r\to 1} d (f(re^{it}), f(0))^p dt
\end{equation*}
is  finite,   which means that $\liminf_{r\to1}d(f(re^{it}),f(0))$ is finite for a.e. $t\in [-\pi,\pi]$. By our assumption on the distance $d$,  we
must have $f_b (t)\in D$ for a.e. $t\in [-\pi,\pi]$.

We   now proceed to show that  $f_b\in\rho\Lambda^p_\alpha$. For $r\in(0,1)$, $t\in [-\pi,\pi]$, and sufficiently small $s\in (0,\frac r2)$, let us
consider the following  piecewise  $\mathcal{C}^1$-smooth curve  contained in the unit disk:
\begin{equation*}
\gamma (x) =
\begin{cases}
(r -x) e^{it}, & 0\le x\le s;\\
(r-s) e^{i\lambda(x)},&  s\le x\le r-s;\\
xe^{i(t+s)}, & r-s\le x\le r,
\end{cases}
\end{equation*}
where
\begin{equation*}
\lambda (x) = \frac {s}{r-2s} (x-s) +t
\end{equation*}
is the affine mapping that maps $[s,r-s]$ onto  $[t,t+s]$. For this  path we have
\begin{equation*}\begin{split}
d (f(r e^{i(t+s)}), f( re^{i t})) &\le
\int _\gamma f^\ast (z)|dz|  = -\int _{0}^s f^\ast ( (r-x)  e^{it }) dx
\\&+ \int _s^{r-s}f^\ast ( (r-s )  e^{i \lambda(x)}) (r-s)\lambda'(x) dx
\\& +\int _{r-s}^r f^\ast (x  e^{i (t+s) }) dx.
\end{split}\end{equation*}
Let $\chi_A$ be the characteristic function of the set $A=\{e^{i\lambda}:t\le\lambda\le t+s\}\subseteq\partial\mathbb {D}$.     The second integral
above may be estimated  applying  the H\"{o}lder inequality in the following  way
\begin{equation*}\begin{split}
\int _s^{r-s} f^\ast  ( (r-s )  e^{i \lambda(x)}) \lambda'(x) dx & =
\int _t^{t+s}  f^\ast  ( (r-s )  e^{i \lambda}) d\lambda
\\& = \int_{-\pi}^{\pi} f^\ast  ( (r-s )  e^{i \lambda})\chi_A(\lambda) d\lambda
\\&\le\left(\int_{-\pi}^{\pi} f^\ast  ( (r-s )e^{i \lambda})^pd\lambda\right)^{\frac 1p}
\left (\int_{-\pi}^{\pi}\chi_A (\lambda)^q \right)^{\frac 1q}
\\& = m_p (r-s, f^\ast_\rho)\cdot s,
\end{split}\end{equation*}
which results
\begin{equation*}\begin{split}
d (f(r e^{i(t+s)}), f( re^{i t})) &\le
\int _{r-s}^r f^\ast (x  e^{it }) dx +\int _{r-s}^r f^\ast (x  e^{i (t+s) }) dx
\\&+ s(r-s)  m_p (r-s, f^\ast).
\end{split}\end{equation*}
Applying now the elementary inequality $(\alpha+\beta+\gamma)^p\le 3^{p-1}(\alpha^p+\beta^p+\gamma^p)$,     where $\alpha$, $\beta$,  $\gamma$  are
non-negative numbers, then integrating against $t$ over $[-\pi,\pi]$,                                           and finally applying the inequality
$(\alpha +\beta+ \gamma)^{\frac 1p}\le \alpha^{\frac 1p}+\beta^{\frac 1p}+\gamma^{\frac 1p}$, we obtain
\begin{equation*}\begin{split}
\left(\int_{-\pi}^{\pi} d (f(r e^{i(t+s)}), f( re^{i t}))^p  {dt}\right)^{\frac 1p}&
\le 3^{1-\frac 1p} \left( \int_{-\pi}^{\pi} \left( \int _{r-s}^r f^\ast (x  e^{it }) dx\right)^p {dt}\right)^{\frac 1p}
\\&+3^{1-\frac 1p}\left(\int_{-\pi}^{\pi} \left( \int _{r-s}^r f^\ast(x  e^{i(t+s) })dx\right)^p {dt} \right)^{\frac 1p}
\\&+3^{1-\frac 1p} s (r-s)  m_p (r-s, f^\ast)\cdot (2\pi)^\frac 1p.
\end{split}\end{equation*}
Applying the Minkowski inequality on the first and the second  integral on the right side of the last inequality,  we obtain
\begin{equation*}\begin{split}
\left(\int_{-\pi}^{\pi} d  (f(r e^{i(t+s)}), f( re^{i t}))^p  {dt}\right)^{\frac 1p} & \le
 6 \int _{r-s}^r m_p(x,f^\ast) {dx}
 + 27 s (r-s) m_p (r-s, f^\ast).
\end{split}\end{equation*}
Now, having in mind  \eqref{EQ.CALPHA}, we derive
\begin{equation*}\begin{split}
\left(\int_{-\pi}^{\pi} d  (f(r e^{i(t+s)}), f( re^{i t}))^p  {dt}\right)^{\frac 1p}
& \le \frac {6C} {\alpha} ( (1- (r-s))^\alpha -(1-r)^\alpha ) \\&+ 27 C s (r-s) (1-(r-s))^{\alpha-1}.
\end{split}\end{equation*}
If we let $r\to 1$, by applying the Fatou lemma on the term on the left side, we obtain
\begin{equation*}\begin{split}
\left(\int_{-\pi}^{\pi} d  (f_b(t+s), f_b(t) )^p  {dt}\right)^{\frac 1p}&
\le \frac {6C} {\alpha} s ^\alpha + 27 C (1-s) s^\alpha\le \frac {33 C} \alpha s^\alpha.
\end{split}\end{equation*}
Therefore,
\begin{equation*}\begin{split}
\sup_{s\in (0,h)} \left(\int_{-\pi}^{\pi} d  (f_b(t+s), f_b(t) )^p  {dt}\right)^{\frac 1p}
\le\frac { 33 C} \alpha h^\alpha,
\end{split}\end{equation*}
which   means that   $f_b\in \rho \Lambda_\alpha^p$.
\end{proof}

\section{Hardy  and  Littlewood theorems and invariant metrics}

\begin{definition}\label{DEF.WEAK}
We say that a domain $D$ is $\mathcal{F}-$transitive, $\mathcal{F}\subseteq\mathrm{Aut}(D)$, if there exists a compact set $K\subseteq D$ such that for
every $\zeta\in D$ there exists $\varphi = \varphi _\zeta \in\mathcal{F}$  with  $\varphi (\zeta)\in K$.
\end{definition}

\begin{definition}
We say that a metric $\rho$ on a domain $D $ is $\mathcal {F}$-invariant, $\mathcal {F}\subseteq\mathrm{Aut}(D)$, if the corresponding distance $d = d_\rho$
satisfies
\begin{equation}\label{EQ.DRHO}
d (\varphi(\zeta), \varphi(\eta)) = d(\zeta,\eta),\quad  \zeta,  \eta \in D
\end{equation}
for every $\varphi \in\mathcal {F}$.
\end{definition}

Note that  if a metric  $\rho$ on a domain $D$  is  $\mathcal{F}$-invariant, then  we have
\begin{equation}\label{EQ.DERIVATIVE}
|\varphi'(\zeta)|  = \frac{\rho  (\zeta)}{\rho  (\varphi(\zeta))},\quad \zeta\in D
\end{equation}
for every  $\varphi \in \mathcal {F}$. Indeed,  for different $\zeta$ and $\eta$ the relation \eqref{EQ.DRHO} may be rewritten as
\begin{equation*}
\frac {d (\varphi(\zeta),\varphi (\eta))} {|\varphi(\zeta) -\varphi(\eta)|}
\left|\frac{ \varphi(\zeta)-\varphi(\eta) } {\zeta-\eta }\right|
\left( \frac {d(\zeta,\eta)}{|\zeta-\eta|} \right)^{-1} = 1.
\end{equation*}
If we let   $\eta \to \zeta$ above, we obtain
\begin{equation*}
\rho  (\varphi (\zeta)) | \varphi'(\zeta)| {\rho (\zeta)} ^{-1}=1,
\end{equation*}
which gives \eqref{EQ.DERIVATIVE}.

The following  result is the converse of  Theorem \ref{TH.2}, but  in a special case.

\begin{theorem}\label{TH.3}
Assume that a domain $D$ is $\mathcal{F}$-transitive, $\mathcal{F}\subseteq \mathrm{Aut} (D)$, and let $\rho$  be  a metric on $D$ which is
$\mathcal {F}$-invariant. If for an analytic mapping $f:\mathbb{D}\to D$ we have $f\in\rho H^p$, $p\in[1,\infty)$,  and $f_b\in\rho \Lambda^p_\alpha$, $\alpha\in(0,1]$,  then
\begin{equation*}
m_p(r, f^\ast)  =  \mathcal {O} (1-r)^{\alpha-1},\quad r\to 1.
\end{equation*}
\end{theorem}

Clearly, every simply--connected domain $D\subseteq\mathbb{C}$ with at least two boundary points is $\mathrm{Aut}(D)$-transitive. If the domain $D$
is bounded, it is well known that the Bergman metric is $\mathrm {Aut}(D)$-invariant. Hence, the preceding theorem gives the converse  part of
Theorem  \ref{TH.BERGMAN}.

The domain $\mathbb{C}$ is $\{\varphi_b(z)=z+b:b\in\mathbb{C}\}$-transitive; for the role of the compact set one may take any  one-point set. Note
that this domain is also  $\{\varphi_{n,m}(z)=z+m+in:m,n\in\mathbb{Z}\}$-transitive; here  once  should take $[-1,1]\times [-1,1]$  for the compact set.
If we take $\rho\equiv 1$ on $ \mathbb{C}$, then $d_\rho$ is  the Euclidean distance which is $\{\varphi_b(z)=z+b:b\in\mathbb{C}\}$-invariant.
Therefore, in this spacial case  Theorem \ref{TH.3} is  the converse  part of the Hardy and Littlewood theorem given in Proposition \ref{PROP.HL2}.

\begin{proof}[Proof of Theorem \ref{TH.3}]
Let $K$ be a  compact set as  in Definition \ref{DEF.WEAK}, and denote
\begin{equation*}
C_ K = \max_{\zeta \in K}  \rho  (\zeta).
\end{equation*}

Let $z\in\mathbb{D}$ and select  a conformal mapping $\varphi = \varphi_z\in\mathcal{F}$ (this mapping depends on $z$) such that $\varphi (f(z))\in K$.
In the sequel we shall consider the analytic mapping $g = g_z = \varphi\circ f:\mathbb {D}\to D$. It is easy to see (by  the $\mathcal {F}$-invariance
of the  distance $d$) that $g\in \rho H^p$ and $g_b\in\rho \Lambda^p_\alpha$.

Using \eqref{EQ.DERIVATIVE} we obtain
\begin{equation*}\begin{split}
f^\ast  (z)  & = {\rho (f(z))}  |f' (z)| = \rho (\varphi(f(z))) \frac{ \rho (f(z))}{ \rho  (\varphi (f(z)) )}  |f' (z)|
\\& =  \rho (\varphi(f(z)))|\varphi'(f(z))||f'(z)|\le C_K |(\varphi\circ f)' (z) |
  = C_K |g' ( z )|.
\end{split}\end{equation*}
Therefore, we have proved
\begin{equation}\label{EQ.FSTAR}
f^\ast (z) \le C_K |g' ( z )|, \quad z\in \mathbb {D}.
\end{equation}

Let $z=re^{i\theta}$, $r\in (0,1)$, $\theta\in [-\pi,\pi]$. We are going to show  that   the $\rho$-derivative of $f$ satisfies the appropriate growth
condition. We separate the following two cases: (i) $\alpha\in (0,1)$; here  we need to prove that  $m_p (r,f^\ast) = \mathcal {O}(1-r)^{\alpha-1}$,
$r\to 1$, and (ii) $\alpha = 1$,  where  we have to show that $m_p (r,f^\ast)$ is bounded in $r\in (0,1)$.

(i) By the Cauchy integral theorem applied to  $g - g_b (\theta)$, from  \eqref{EQ.FSTAR} we obtain
\begin{equation*}\begin{split}
& f^\ast  (re^{i\theta})\le C_ K g'(re^{i\theta})
\\& =  C_K \left|\int _{-\pi}^{\pi} \frac{g_b(t)-g_b(\theta)}{(e^{it} -re^{i\theta})^2}e^{it} \frac {dt} {2\pi}\right|
\le C_K  \int _{-\pi}^{\pi} \frac{|g_b(t) - g_b(\theta)|}{|e^{it} -re^{i\theta}  |^2} \frac {dt} {2\pi  }
\\&\le C_K  C_\rho^{-1} \int _{-\pi}^{\pi} \frac{d (g_b(t), g_b(\theta) ) }{ 1-2r \cos(t-\theta) +r^2}\frac {dt} {2\pi  }
= \frac{C_K  C_\rho^{-1}}{2\pi}  \int _{-\pi}^{\pi}\frac{d(g_b(t+\theta) , g_b(\theta) ) dt}{ 1-2r \cos t +r^2 }
\end{split}\end{equation*}
for a.e. $\theta\in [-\pi,\pi]$. Since the distance $d$ is invariant with respect to the conformal  mapping $\varphi$, for $s\in (0,1)$ we have
\begin{equation*}
d (g (s e^{i(t+\theta)}), g(s e^{i\theta})) =  d (\varphi (f(se^{i(t+\theta)})),\varphi (f (s e^{i\theta}) ) )
= d (f(se^{i(t+\theta)}), f(se^{i{\theta}})).
\end{equation*}
Taking   the limit as $s\to 1$, it follows
\begin{equation*}
d (g_b (t+\theta), g_b( \theta ) ) = d (f_b(t+\theta), f_b(\theta ))\  \text{for a.e.}\ t\in [-\pi,\pi].
\end{equation*}
Therefore, we  may conclude that the following inequality holds
\begin{equation*}
f^\ast(r e^{i \theta})\le \frac{C_K  C_\rho^{-1}}{2\pi} \int_{-\pi}^{\pi}\frac{d (f_b(t+\theta),f_b(\theta))} {1-2r\cos t+r^2} {dt}.
\end{equation*}
Now, applying   the    Minkowski inequality,  we obtain
\begin{equation*}\begin{split}
\left(\int_{-\pi}^{\pi}f^\ast (r e^{i \theta}) ^p {d\theta} \right)^{\frac  1p}
&\le \frac{C_K  C_\rho^{-1}}{2\pi}
\left(\int_{-\pi}^{\pi}\left( \int _{-\pi}^{\pi}
\frac{d  (f_b( t+\theta) , f_b(\theta) ) }{1-2r\cos t+r^2}{dt} \right)^p {d\theta} \right)^{\frac  1p}
\\&\le \frac{C_K  C_\rho^{-1}}{2\pi}
\int_{-\pi}^{\pi}\left(\int _{-\pi}^{\pi}
\left(\frac{d (f_b(t+\theta), f_b(\theta) ) } {1-2r\cos t+r^2}\right)^p  {d\theta} \right)^{\frac  1p}  {dt}.
\end{split}\end{equation*}
Since $f_b\in \rho\Lambda^p_\alpha$, there exist a constant $C_\alpha$ such that
\begin{equation*}
\left(\int _{-\pi}^{\pi}
d(f_b(\theta+t) , f(\theta) )^p  {d\theta} \right)^{\frac  1p}  {dt}\le C_\alpha |t|^\alpha.
\end{equation*}
Having in mind  the  inequality
\begin{equation*}
1-2r\cos t +r^2 = (1-r)^2+4r\sin ^2 \frac t2 \ge (1-r)^2  + \frac {4rt^2}{\pi^2},
\end{equation*}
it follows
\begin{equation*}\begin{split}
\left(\int_{-\pi}^{\pi} f^\ast (r e^{i \theta}) ^p   {d\theta} \right)^{\frac  1p} &
\le \frac{C_K  C_\rho^{-1} C_\alpha}{2\pi} \int_{-\pi}^{\pi} \frac {|t|^\alpha}{1-2r\cos t+r^2} {dt}
\\& = \frac{C_K  C_\rho^{-1} C_\alpha}{\pi}  \int_{0}^{\pi}\frac {t^\alpha dt}{1-2r\cos t+r^2}
\\& \le\frac{C_K  C_\rho^{-1} C_\alpha}{\pi}
\int_{0}^{\pi}\frac {t^\alpha dt}{(1-r)^2 +\frac{4rt^2}{\pi^2} }\quad (u=\frac {2\sqrt{r}}{\pi} \frac{t}{1-r})
\\& =\frac{C_K  C_\rho^{-1} C_\alpha}{\pi}
\left(\frac {\pi}{2\sqrt{r}}\right)^{\alpha+1}(1-r)^{\alpha-1}\int_{0}^{\frac {2\sqrt{r}}{1-r}}\frac {u^\alpha}{1+u^2}du
\\& \le\frac{C_K  C_\rho^{-1} C_\alpha}{\pi}
\left(\frac {\pi}{2\sqrt{r}}\right)^{\alpha+1}(1-r)^{\alpha-1}\int_{0}^{\infty}\frac {u ^\alpha dt}{1+ u^2}du
\\& = \frac{C_K  C_\rho^{-1} C_\alpha}{\pi}
\left(\frac {\pi}{2\sqrt{r}}\right)^{\alpha+1} (1-r)^{\alpha-1} \frac \pi 2 \frac1 {\cos(\frac {\alpha\pi}2)},
\end{split}\end{equation*}
which  proves that   $m_p (r, f^\ast) = \mathcal {O}(1-r)^{\alpha-1}$, $r\to 1$.

(ii) Here we shall first note that $f_b$ is absolutely continuous. Since $f_b\in\rho\Lambda_1^p\subseteq \rho\Lambda_1^1\subseteq \Lambda_1^1$, by
the  Hardy and Littlewood  result \cite[Lemma 1, p. 72]{DUREN.BOOK} it follows that $f_b$ coincides a.e. with a function of bounded variation on
$[-\pi,\pi]$. Now, by \cite[Theorem 3.11 and Theorem 3.10]{DUREN.BOOK} it follows that $f_b$ is absolutely continuous  and that $f'$ belongs to
$H^1$. Moreover, by the same theorems,  we have the following  connection  between the derivatives
\begin{equation*}
f_b'(t) = i e^{it} f'(e^{it})\ \text {for a.e.}\ t\in [-\pi,\pi].
\end{equation*}
In the sam way,  because   $g\in \rho H^p$ and $g_b \in \rho \Lambda_\alpha^p$, we conclude that $g_b$ is also absolutely continuous and $g'\in H^1$.

Since  $f_b$ and $g_b$ are absolutely  continuous, and $f_b\in D$ for a.e.  $t\in[-\pi,\pi]$, we may conclude that $g_b(t)=(\varphi \circ f_b)(t)$
and $g _b '(t)=\varphi '(f_b( t ))  f_b'(t)$ for  almost every $t\in [-\pi,\pi]$.

Since $f_b$ belongs to the  class  $\rho\Lambda_1^p$, there exists a constant $C$ such that
\begin{equation*}
\int_{-\pi}^{\pi}\frac {d(f_b(t+h),f_b(t) )^p }{| h |^p}dt\le C
\end{equation*}
for  $h $  sufficiently close to  $0$. If we let $h\to 0$ above,  by the Fatou theorem we obtain
\begin{equation}\label{EQ.BETA.INT}
\int_{-\pi}^{\pi} f_b^\ast (t) ^p dt\le C,
\end{equation}
i.e., $f_b^\ast\in L^p [-\pi,\pi]$. For a.e.  $t\in [-\pi,\pi]$ we have
\begin{equation*}\begin{split}
|g_b'(t)| & = |\varphi'(f _b (t))|  |f_b'(t)| =\frac{\rho (f _b (t)) } {\rho (\varphi (f_b(t)))}  \frac {f_b^\ast (t)}{ \rho (f_b (t))}
=\frac{ {f_b^\ast (t)} } {\rho  (\varphi (f_b(t)))} \le C _\rho^{-1} f_b^\ast (t),
\end{split}\end{equation*}
since   $\rho  \ge C_\rho$ on $D$.

We have $g ' (w) = \mathrm {P} [g'(e^{it})] (w) = \mathrm {P} [-ie^{-it}g_b'(t)](w)$ (recall that $g\in H^1$). Applying now  the  Jensen
inequality and the  inequality above:  $|g_b'(t)|\le C _\rho^{-1} f_b^\ast (t) $, we obtain
\begin{equation*}\begin{split}
|g'(w)|^p & \le  |\mathrm {P}[-i e^{-it} g_b'(t)](w)|^p
\le \mathrm {P}[|g_b'|^p]|(w)\le C_\rho^{-p} \mathrm {P}[{f_b^\ast} ^p]| (w),\quad w\in \mathbb{D}.
\end{split}\end{equation*}
Particularly, for  $w = z$, having in mind the inequality \eqref{EQ.FSTAR},  we derive
\begin{equation*}
f^\ast (z)^p   \le C_K^p   | g'(z)|^p \le C_K^p  C_\rho^{-p} \mathrm {P}[{f_b^\ast}(t)^p] (z).
\end{equation*}
Applying the Fubini theorem we obtain
\begin{equation*}
\int_{-\pi}^{\pi} f^\ast (re^{i\theta})^p {d\theta}
\le  C_K^p C_\rho^{-p} \frac 1{2\pi} \int_{-\pi}^{\pi} \left(\int_{-\pi}^{\pi}
\frac {1-r^2}{1-2r\cos (t-\theta)  +r^2} {d\theta}  \right) f_b^{\ast}(t)^pdt.
\end{equation*}
Since the  inner integral is equal to $2\pi$, the preceding inequality may be rewritten as
\begin{equation*}
m_p (r,f^\ast)\le  C_K C_\rho^{-1}   \|f_b^\ast\|_p,
\end{equation*}
which implies that  $m_p(r, f^\ast)$ is bounded  in $r \in (0,1)$.
\end{proof}

\section*{Acknowledgements}
I would like to thank the referee of this work for carefully reading and pointing out several errors in previous versions of the manuscript.


\begin{thebibliography}{10}

\bibitem{BLOOM.CA}
S. Bloom  and G.S. De Souza, \textit{Weighted Lipschitz Spaces and Their Analytic  Characterizations},
Constr. Approx. \textbf{10} (1994) 339--376.

\bibitem{DUREN.BOOK}
P.L. Duren, \textit{Theory of HP spaces},
Academic Press, New York and  London, 1970.

\bibitem{HL.MZ}
G.H. Hardy and J.E. Littlewood, \textit{A convergence criterion for Fourier series},
Math. Z. \textbf{28} (1928), 612--634.

\bibitem{KRANTZ.BOOK}
St.G. Krantz, \textit{Complex Analysis: The Geometric Viewpoint},
MAA,  Washington, 2004.

\bibitem{MARKOVIC.JGA}
M. Markovi\'{c}, \textit{Representations for the Bloch Type Semi-norm of Fr\'{e}chet Differentiable Mappings},
J. Geom. Anal. \textbf{31} (2021), 7947--7967.

\bibitem{ZHU.JLMS}
K. Zhu, \textit{Distances and Banach spaces of holomorphic functions on complex domains},
J. Lond. Math. Soc. \textbf{49} (1994), 163--182.

\bibitem{NIKOLOV.COLLM}
N. Nikolov and  M. Trybula, \textit{Estimates of the Bergman distance on Dini--smooth bounded planar domains},
Coll. Math. \textbf{67} (2016), 407--414.

\bibitem{PAVLOVIC.BOOK}
M. Pavlovi\'{c}, \textit{Function classes on the unit disc},
De Gruyter Studies in Mathematics, 2019.

\bibitem{YAMASHITA.TAMS}
Sh. Yamashita, \textit{Smoothness of the boundary values of functions bounded and holomorphic in the disk},
Trans. Amer. Math. Soc. \textbf{272} (1982), 539--544.

\end{thebibliography}
\end{document}